\documentclass[oneside,11pt]{amsart}

\usepackage{float}
\usepackage{amscd,amssymb}
\usepackage[mathscr]{eucal}
\usepackage[all]{xy}
\usepackage{mathrsfs}
\usepackage{amsfonts}
\usepackage{amsmath}
\usepackage{amsthm}
\usepackage{latexsym}

\title{A vanishing theorem for differential operators in positive
  characteristic}

\author{Alexander Samokhin}

\address{Institute for Information Transmission Problems, Moscow,
  Russia}

\email{alexander.samokhin@gmail.com}

\newcommand{\Oo}{\mathcal O}
\newcommand{\Uu}{\mathcal U}

\newcommand{\Pp}{\mathbb P}

\newcommand{\Ff}{\mathcal F}
\newcommand{\Ee}{\mathcal E}
\newcommand{\Ll}{\mathcal L}

\newcommand{\D}{\mathcal D}

\newcommand{\T}{\mathcal T}

\newcommand*{\RHom}{\mathop{\mathrm RHom}\nolimits}
\newcommand*{\Hom}{\mathop{\mathrm Hom}\nolimits}
\newcommand*{\Dd}{\mathop{\mathrm D\kern0pt}\nolimits}
\newcommand*{\DD}{\mathop{\mathbb D\kern0pt}\nolimits}

\newcommand*{\Ext}{\mathop{\mathrm Ext}\nolimits}

\newtheorem{theorem}{Theorem}[section]
\newtheorem{corollary}{Corollary}[section]
\newtheorem{lemma}{Lemma}[section]

\newtheorem{remark}{Remark}[section]
\newtheorem{definition}{Definition}[section]

\numberwithin{equation}{section}

\long\def\comment#1{}

\begin{document}

\maketitle

\begin{abstract}
Let $X$ be a smooth variety over an algebraically closed field $k$ of
characteristic $p$, and ${\sf F}\colon X\rightarrow X$ the Frobenius
morphism. We prove that if $X$ is an incidence variety (a partial flag
variety in type ${\bf A}_n$) or a smooth quadric (in this case $p$ is
supposed to be odd) then ${\rm H}^{i}(X,{\mathcal End}({\sf F}_{\ast}\Oo _X)) = 0$ for $i
> 0$. Using this vanishing result and the derived localization theorem for crystalline
differential operators \cite{BMR}, we show that the Frobenius direct
image ${\sf F}_{\ast} \Oo _X$ is a tilting bundle on these varieties provided that $p>h$, the Coxeter number of the corresponding group.
\end{abstract}


\section*{Introduction}
Let $X$ be a smooth variety over an algebraically closed field $k$ of
arbitrary characteristic, and $\D _X$ the sheaf of differential operators \cite{Gro}. Recall
that the sheaf $\D _X$ is equipped with a filtration by degree of an
operator such that the associated graded sheaf is isomorphic to
$\bigoplus ({\sf S}^{i}\Omega ^1_X)^{\ast}$. Our goal 
is to study the cohomology vanishing of
sheaves of differential operators in positive 
characteristic. We are primarily interested in the case when $X$ is a
homogeneous space of a semisimple algebraic group ${\bf G}$ over $k$, i.e. $X =
{\bf G}/{\bf P}$ for some parabolic subgroup ${\bf P}\subset {\bf
  G}$. Below we survey some vanishing and non-vanishing results for
sheaves of differential operators. If $X$ is a smooth variety over a
field $k$ of characteristic zero then there is an
isomorphism $({\sf S}^i\Omega ^1_X)^{\ast}={\sf S}^{i}\T _X$ for
$i\geq 0$, and it is known that symmetric powers of the tangent bundle of ${\bf G}/{\bf P}$ have vanishing higher cohomology groups \cite{BK}. In particular,
this implies the vanishing of higher cohomology of the sheaf $\D _{{\bf
  G}/{\bf P}}$. This vanishing is also a consequence of the
Beilinson--Bernstein localization theorem \cite{BB}. If the field
$k$ is of characteristic $p>0$ then the situation is less clear. A
crucial tool to prove the vanishing theorem for symmetric powers of
the tangent bundle of ${\bf G}/{\bf P}$ in characteristic zero is the Grauert--Riemenschneider
theorem, which does not hold in characteristic $p$ in general. 
The characteristic $p$ counterpart of the cohomology vanishing of
symmetric powers of the tangent bundle of ${\bf G}/{\bf P}$ is inaccessible at the moment, though it is
believed to be true (e.g., \cite{BrKum}). There are cases when such a theorem is known -- for example, for flag
varieties ${\bf G}/{\bf B}$ in good characteristic \cite{KLT}. These
results, however, do not give any information about cohomology groups of the
sheaf $\D _{{\bf G}/{\bf P}}$, since for a variety $X$ in positive 
characteristic the sheaves of graded rings
$\bigoplus ({\sf S}^{i}\Omega ^1_X)^{\ast}$ and $\bigoplus
{\sf S}^{i}\T _X$ are no longer isomorphic. The latter sheaf is
isomorphic to the associated graded ring of the sheaf of crystalline
differential operators ${\rm D}_{X}$ \cite{BMR}, and the vanishing theorem in
the case of flag varieties gives that ${\rm H}^{i}({\bf G}/{\bf
  B},{\rm D}_{{\bf G}/{\bf B}}) = 0$ for $i>0$ and good $p$. On the
other hand, some non-vanishing results are known for both the sheaves
$\D _{{\bf G}/{\bf B}}$ and $\bigoplus ({\sf S}^{i}\Omega
^1_{{\bf G}/{\bf B}})^{\ast}$. It follows from \cite{KaLau} that for
the flag variety ${\bf SL}_5/{\bf B}$ the former sheaf has a
non-vanishing higher cohomology group. Moreover, the latter sheaf has
always a non-vanishing cohomology group unless ${\bf G}/{\bf B}$ is
isomorphic to $\Pp ^1$ (\cite{Haas}, \cite{LM}).

From now on let $k$ be an algebraically closed field of characteristic
$p>0$. Recall that for a variety $X$ over $k$ the sheaf of
differential operators $\D _X$ is a union of sheaves of matrix
algebras \cite{Haas}. Precisely, $\D _X = \bigcup _{m\geq 1} {\mathcal
  End}_{\Oo _X^{p^m}}(\Oo _X)$. The main results of this paper are:

\begin{theorem}\label{th:1stth}
Let $X$ be a partial flag variety of type $(1,n,n+1)$ of the group
${\bf SL}_n$. In other words, $X = \{(l\subset H\subset V)\}$, where
$l$ and $H$ are a line and a hyperplane in the vector space $V$ of
dimension $n+1$, respectively. Then ${\rm H}^{i}(X,{\mathcal End}_{\Oo _X^{p}}(\Oo
_X))=0$ for $i>0$.
\end{theorem}

\begin{theorem}\label{th:2ndth}
Let ${\sf Q}_n\subset \Pp ^{n+1}$ be a smooth quadric of dimension
$n$. Assume that $p$ is an odd prime. Then ${\rm H}^{i}({\sf Q}_n,{\mathcal End}_{{\Oo ^{p}_{{\sf Q}_n}}}(\Oo _{{\sf Q}_n}))=0$
for $i>0$.
\end{theorem}

For a variety $X$ the sheaf ${\mathcal End}_{\Oo _X^{p}}(\Oo _X)$ is
also called the ``sheaf of small differential operators'' (as it is isomorphic to
the central reduction of the sheaf ${\rm D}_X$), which explains the
title of the paper. The proof of Theorems \ref{th:1stth} and \ref{th:2ndth} uses properties of sheaves of crystalline differential operators
\cite{BMR} and vanishing theorems for line bundles on the cotangent
bundles of ${\bf G}/{\bf P}$. Independently, in a recent paper
\cite{Lan} A. Langer, using different methods, described the decomposition of Frobenius push-forwards ${{\sf
F}_m}_{\ast}\Oo _{{\sf Q}_n}$ into a direct sum of indecomposable
bundles. This gives another proof of Theorem \ref{th:2ndth}.

Our \ interest \ in \ the \ vanishing \ of \ higher \ cohomology \
  groups \ of  \ sheaves ${\mathcal
  End}_{\Oo ^p_{{\bf G}/{\bf P}}}(\Oo _{{\bf G}/{\bf P}})$ is twofold.
  First, as we saw above, it is related to the vanishing of higher cohomology of the
  sheaf $\D _{{\bf G}/{\bf P}}$. In the case when ${\bf P}={\bf B}$,
  the Borel subgroup, the vanishing of higher cohomology groups of $\D
  _{{\bf G}/{\bf B}}$ implies ${\sf D}$-affinity of the flag variety
  ${{\bf G}/{\bf B}}$ (\cite{AnKan}, \cite{Haas}, \cite{Lan},
  \cite{SamD-aff}). On the other hand, it has implications for the
  derived category of coherent sheaves: if the
  vanishing holds for a given ${\bf G}/{\bf P}$ then the derived
  localization theorem \cite{BMR} implies that the Frobenius
  pushforward ${\sf F}_{\ast}\Oo _{{\bf G}/{\bf P}}$ is a tilting
  bundle on ${\bf G}/{\bf P}$, provided that $p>h$, the Coxeter number
  of $\bf G$. We discuss these questions in a greater detail in
  Section \ref{sec:derequiv} (see also \cite{HKR}, \cite{KY}, \cite{Lan}, \cite{Samtiltfrob}).

We would like to conclude the introductory part with a conjecture. There are a number of previously known results (see the list
below), to which the present paper gives an additional evidence,
that suggest that one could hope for the vanishing of higher cohomology
of the sheaves of small differential operators on homogeneous spaces
${\bf G}/{\bf P}$.

{\bf Conjecture}: Let ${\bf G}/{\bf P}$ be the homogeneous
space of a semisimple simply connected algebraic group $\bf G$ over
$k$. Then \ for \ sufficiently \ large \ $p$ \ one \ has \ ${\rm H}^{i}({\bf G}/{\bf P},{\mathcal
  End}_{\Oo _{{\bf G}/{\bf P}}^{p}}(\Oo _{{\bf G}/{\bf P}}))=0$ for
$i>0$.

Surprisingly enough, this cohomology vanishing is known, apart from
the above theorems, in very few cases, namely for projective spaces, the flag variety ${\bf
  SL}_3/{\bf B}$ by \cite{Haas}, the flag variety in type ${\bf B}_2$
by \cite{AnKan}. Let us also remark that the counterexample to the ${\sf D}$-affinity from 
\cite{KaLau}, which is based on the study of a particular ${\mathcal D}$-module on
the grassmannian ${\rm Gr}_{2,5}$, does not rule out the possibility
of the conjectural vanishing in this case. The first step
to further test the conjecture is to check it for ${\rm Gr}_{2,5}$;
this will be a subject of a forthcoming paper \cite{SamFrobhom}.

\subsubsection*{Acknowledgements}
We are much indebted to R. Bezrukavnikov for many valuable conversations, his patient explanations of the results of
\cite{BMR}, and sharing his ideas. The present paper was directly
inspired by \cite{BMR} and \cite{BMRII}. We would also like to thank 
M. Finkelberg for the interest and pointing out a gap in an earlier version of the
paper, A. Kuznetsov for his critical remarks and generous help,
B. Keller for his support, and A. Langer for many useful
conversations. An early draft of this paper was written in the fall of
2006 and first appeared in \cite{SamokhinIHES}. The present paper was completed during the
author's stays at the Institute of Mathematics of Jussieu, the IHES, and the Max-Planck Institute for Mathematics in Bonn. It is a
great pleasure to thank all these institutions for the hospitality and
support. 


\renewcommand\thesubsection{\arabic{subsection}}

\section{Preliminaries}\label{sec:prelim}
\subsection{The Frobenius morphism}\label{subsec:Frobmorph}

Throughout we fix an algebraically closed field $k$ of characteristic
$p>0$. Let $X$ be a scheme over $k$. The absolute Frobenius morphism ${\sf F}_X$ is an endomorphism
of $X$ that acts identically on the topological space of $X$ and
raises functions on $X$ to the $p$-th power:
\begin{equation}\label{eq:defofFrobmor}
{\sf F}_X\colon X\rightarrow X,\qquad f\in \Oo _X\mapsto f^p\in \Oo _X.
\end{equation}

Denote $X'$ the scheme obtained by base change from $X$
along the Frobenius morphism ${\sf F}\colon {\rm Spec}(k)\rightarrow {\rm
Spec}(k)$. Then the relative Frobenius morphism ${\sf F}\colon X\rightarrow
X'$ is a morphism of $k$-schemes. The schemes $X$ and $X^{'}$ are
isomorphic as abstract schemes. 

For a quasicoherent sheaf $\Ee$ on $X$ one has, the Frobenius
morphism being affine:
$$
{\rm H}^{i}(X,\Ee) = {\rm H}^{i}(X',{\sf F}_{\ast}\Ee).
$$

\subsection{Koszul resolutions}\label{subsec:Kosres}

Let $\sf V$ be a finite dimensional vector
space over $k$ with a basis $\{e_1,\dots ,e_n\}$. Recall that the {\it $r$-th exterior power} $\bigwedge
^{r} \sf V$ of $\sf V$ is defined to be the $r$-th tensor power ${\sf V}^{\otimes
  r}$ of $\sf V$ divided by the vector subspace spanned by the elements:
$$ u_1\otimes \dots \otimes u_r - (-1)^{{\rm sgn} \sigma}u_{\sigma
  _1}\otimes \dots \otimes u_{\sigma (r)}$$

for all the permutations $\sigma \in \Sigma _r$ and
  $u_1,\dots ,u_r \in V$. Similarly, the {\it $r$-th symmetric power}
  ${\sf S}^{r}\sf V$ of $\sf V$ is defined to be the $r$-th tensor power
  ${\sf V}^{\otimes r}$ of $\sf V$ divided by the vector subspace spanned by the
  elements 
$$ u_1\otimes \dots \otimes u_r - u_{\sigma
  _1}\otimes \dots \otimes u_{\sigma (r)}$$

for all the permutations $\sigma \in \Sigma _r$ and
  $u_1,\dots ,u_r \in \sf V$.\\

Let 
\begin{equation}\label{eq:seqvecsp}
0\rightarrow {\sf V}'\rightarrow {\sf V}\rightarrow {\sf V}''\rightarrow 0
\end{equation}

be a short exact sequence of vector spaces. For any $n > 0$
there is a functorial exact sequence (the Koszul resolution,
(\cite{Jan}, II.12.12))

\begin{equation}\label{eq:Kosqresol}
\dots \rightarrow {\sf S}^{n-i}{\sf V}\otimes \bigwedge ^i{\sf V}'\rightarrow \dots
\rightarrow {\sf S}^{n-1}{\sf V}\otimes {\sf V}'\rightarrow {\sf
  S}^n{\sf V}\rightarrow
{\sf S}^n{\sf V}''\rightarrow 0.
\end{equation}

Another fact about symmetric and exterior powers is the following (\cite{Har}, Exercise
5.16). For a short exact sequence such as (\ref{eq:seqvecsp}) one has for each $n$ the
filtrations

\begin{equation}
{\sf S}^{n}{\sf V} = {\rm F}_n\supset {\rm F}_{n-1}\supset \dots \quad \mbox{and} \quad \bigwedge
^{n}{\sf V} = {\rm F}'_n\supset {\rm F}'_{n-1}\supset \dots
\end{equation}

such that 

\begin{equation}
{\rm F}_i/{\rm F}_{i-1}\simeq {\sf S}^{n-i}{\sf V}'\otimes {\sf
  S}^{i}{\sf V}''
\end{equation}

and 
\begin{equation}
{\rm F}'_i/{\rm F}'_{i-1}\simeq \bigwedge ^{n-i}{\sf V}'\otimes
\bigwedge ^{i}{\sf V}''
\end{equation}

When either ${\sf V}'$ or ${\sf V}''$ is a one-dimensional vector space, these
filtrations on exterior powers of $\sf V$ degenerate into short exact sequences. 
If ${\sf V}''$ is one-dimensional then one obtains:

\begin{equation}
0 \rightarrow \bigwedge^{r}{\sf V}'\rightarrow \bigwedge^{r}{\sf V} \rightarrow \bigwedge
^{r-1}{\sf V}'\otimes {\sf V}''\rightarrow 0.
\end{equation}

Similarly, if ${\sf V}'$ is one-dimensional, the above filtration 
degenerates to give a short exact sequence:

\begin{equation}
0 \rightarrow \bigwedge^{r-1}{\sf V}''\otimes {\sf V}'\rightarrow \bigwedge
^{r}{\sf V} \rightarrow \bigwedge
^{r}{\sf V}''\rightarrow 0.
\end{equation}

\subsection{\bf Differential operators}\label{sec:diffoper}

\subsubsection{``True'' differential operators}
The material here is taken from \cite{Gro} and \cite{Haas}.
Let $X$ be a smooth scheme over $k$. Consider the product $X\times X$ and the diagonal $\Delta \subset
X\times X$. Let ${\mathcal J}_{\Delta}$ be the sheaf of ideals of
$\Delta$. 
\begin{definition}
An element $\phi \in {\mathcal End}_{k}(\Oo _X)$ is called a
differential operator if there exists some integer $n\geq 0$ such that
\begin{equation}
{\mathcal J}_{\Delta}^{n}\cdot \phi = 0.
\end{equation}
\end{definition}

One obtains a sheaf $\D _X$, the sheaf of differential
operators on $X$. Denote ${\mathcal
  J}_{\Delta}^{(n)}$ the sheaf of ideals generated by elements 
$a^n$, where $a\in {\mathcal J}_{\Delta}$. There is a filtration
on the sheaf $\D _X$ given by
\begin{equation}
\D _X^{(n)} = \{\phi \in {\mathcal End}_{k}(\Oo _X)\colon {\mathcal
  J}_{\Delta}^{(n)}\cdot \phi = 0\}.
\end{equation}

Since $k$ has characteristic $p$, one checks:
\begin{equation}\label{eq:p-filtration}
\D _X^{(p^{n})} = {\mathcal End}_{\Oo _X^{p^n}}(\Oo _X).
\end{equation}

Indeed, the sheaf ${\mathcal J}_{\Delta}$ is
generated by elements $a\otimes 1 - 1\otimes a$, where $a\in \Oo _X$, hence the sheaf ${\mathcal J}_{\Delta}^{(p^{n})}$ is
  generated by elements $a^{p^{n}}\otimes 1 - 1\otimes
  a^{p^{n}}$. This implies (\ref{eq:p-filtration}). One also checks that this
  filtration exhausts the whole $\D _X$, so one has (Theorem 1.2.4, \cite{Haas}):
\begin{equation}\label{eq:p-filtration}
\D _X = \bigcup _{n\geq 1}  {\mathcal End}_{\Oo _X^{p^n}}(\Oo _X).
\end{equation}

The filtration from (\ref{eq:p-filtration}) was called the
$p$--filtration in {\it loc.cit}. By definition of the Frobenius morphism
one has ${\rm H}^{i}(X,{\mathcal End}_{\Oo _X^{p}}(\Oo _X)) = {\rm
  H}^{i}(X',{\mathcal End}({\sf F}_{\ast}\Oo _X))$. The sheaf $\D _X$
contains divided powers of vector fields (hence the name ``true'') as
opposed to the sheaf of PD-differential operators $\Dd _X$ that is discussed
in the next paragraph.

\subsubsection{Crystalline differential operators}\label{subsec:Berthdiffoper}

The material of this subsection is taken from \cite{BMR}.
We recall, following {\it loc.cit.} the basic properties of crystalline differential operators (differential
operators without divided powers, or PD-differential
operators in the terminology of Berthelot and Ogus).\\

Let $X$ be a smooth variety, ${\T}^{\ast}_X$ the cotangent bundle, and 
${\rm T}^{\ast}(X)$ the total space of ${\T}^{\ast}_{X}$. Denote $\pi
: {\rm T}^{\ast}(X)\rightarrow X$ the projection.

\begin{definition}
The sheaf $\Dd _X$ of {\it crystalline differential operators}
on $X$ is defined as the enveloping algebra of the
tangent Lie algebroid, i.e.,  for an affine open $U\subset X$ the
algebra $\Dd (U)$ contains the subalgebra $\Oo $ of functions, has an
$\Oo$-submodule identified with the Lie algebra of vector fields
$Vect(U)$ on $U$, and these subspaces generate $\Dd (U)$ subject to
relations $\xi_1\xi_2-\xi_2\xi_1=[\xi_1,\xi_2]\in Vect(U)$ for
$\xi_1, \xi_2\in Vect(U)$, and  $\xi \cdot f -f\cdot \xi =\xi(f)$
for $\xi\in Vect(U)$ and $f\in \Oo (U)$.
\end{definition}

Let us list the basic properties of the sheaf $\Dd _X$ \cite{BMR}:

\begin{itemize}

\item The sheaf of non-commutative algebras ${\sf F}_{\ast}\Dd _X$ has
  a center, which is isomorphic to $\Oo _{{\rm T}^{\ast}(X')}$, the
  sheaf of functions on the cotangent bundle to the Frobenius twist of
  $X$. The sheaf ${\sf F}_{\ast}\Dd _X$ is finite over its center.

\item This makes ${\sf F}_{\ast}\Dd _X$ a coherent sheaf on ${\rm
    T}^{\ast}(X')$. Thus, there exists a sheaf of algebras $\DD
 _X$ on ${\rm T}^{\ast}(X')$ such that $\pi _{\ast}\DD _X = {\sf
 F}_{\ast}{\rm D} _X$  (by abuse of notation we denote the projection ${\rm T}^{\ast}(X')\rightarrow X'$
 by the same letter $\pi$). The sheaf $\DD _X$ is an Azumaya algebra
    over ${\rm T}^{\ast}(X')$ of rank  $p^{2{\rm dim}(X)}$. 

\item There is a filtration on the sheaf ${\sf F}_{\ast}{\rm D}_X$
  such that the associated graded ring ${\rm gr}({\sf F}_{\ast}{\rm
  D}_X)$ is isomorphic to ${\sf F}_{\ast}\pi _{\ast}\Oo
 _{{\rm T}^{\ast}(X)} = {\sf F}_{\ast}{\sf S}^{\bullet}\T _X$.

\item Let $i\colon X'\hookrightarrow {\rm T}^{\ast}(X')$ be the zero section embedding.
Then $i^{\ast}\DD _X$ splits
 as an Azumaya algebra, the splitting bundle being ${\sf F}_{\ast}\Oo
 _X$. In other words, $i^{\ast}\DD _X = {\mathcal End}({\sf
 F}_{\ast}\Oo _X)$.        

\end{itemize}

Finally, recall that the sheaf $\Dd _X$ acts on $\Oo _X$ and that this action is not
faithful in general. It gives rise to a map $\Dd _X\rightarrow \D _X$;
its image is the sheaf of ``small differential operators'' ${\mathcal End}_{\Oo _X^p}(\Oo _X)$.

\subsubsection{Cartier descent}\label{subsec:Cartier}

Let $\Ee$ be an $\Oo _X$-module equipped with an integrable connection
$\nabla$, and let $U\subset X$ be an open subset with local
coordinates $t_1,\dots ,t_d$. Let $\partial _1,\dots ,\partial _d$ be
the local basis of derivations that is dual to the basis $(dt_i)$ of
$\Omega ^1_X$. The connection $\nabla$ is said to have zero
$p$-curvature over $U$ if and only if for any local section $s$ of
$\Ee$ and any $i$ one has $\partial _i^ps=0$. For any $\Oo
_{X'}$-module $\Ee$ the Frobenius pullback ${\sf F}^{\ast}\Ee$ is
equipped with a canonical integrable connection with zero
$p$-curvature. The Cartier descent theorem (Theorem 5.1, \cite{Katz})
states that the functor ${\sf F}^{\ast}$ induces an equivalence
between the category of $\Oo _{X'}$-modules and that of $\Oo
_X$-modules equipped with an integrable connection with zero
$p$-curvature.

On the other hand, if an $\Oo _X$-module $\Ee$ is equipped with  an integrable connection with zero
$p$-curvature then it has a structure of left $\Dd _X$-module. Given
the Cartier descent theorem, we see that the Frobenius pullback ${\sf
  F}^{\ast}\Ee$ of a coherent sheaf $\Ee$ on $X'$ is a left $\Dd _X$-module. We will use the Cartier descent in Section \ref{sec:derequiv}.

\section{\bf Vanishing theorems for line bundles}

Of crucial importance for us are the vanishing theorems for cotangent
bundles of homogeneous spaces. Let $\bf G$ be a connected, simply connected, semisimple algebraic
group over $k$, $\bf B$ a Borel subgroup of $\bf G$, and $\bf T$ a maximal
torus. Let $R({\bf T},{\bf G})$ be the root system of $\bf G$ with
respect to $\bf T$, $R^{+}$ the subset of positive roots, $S\subset
R^{+}$ the simple roots, and $h$ the Coxeter number of $\bf G$ that is
equal to $\sum m_i$, where $m_i$ are the coefficients of the highest
root of $\bf G$ written in terms of the simple roots $\alpha _i$.
 By $\langle \cdot, \cdot \rangle$ we denote the natural pairing $X({\bf
  T})\times Y({\bf T})\rightarrow {\mathbb Z}$, where $X({\bf T})$ is
the group of characters (also identified with the weight lattice) and
$Y({\bf T})$ the group of one parameter subgroups of $\bf T$ (also
identified with the coroot lattice). For a
subset $I\subset S$ let ${\bf P} = {\bf P}_{I}$ denote the associated
parabolic subgroup. Recall that the group of characters $X({\bf P})$
of $\bf P$ can be identified with $\{\lambda \in X({\bf T})|\langle
\lambda ,\alpha ^{\vee}\rangle = 0,$ for all $\alpha \in I\}$. In
particular, $X({\bf B}) = X({\bf T})$. A weight $\lambda \in X({\bf
  B})$ is called {\it dominant} if $\langle \lambda ,\alpha
^{\vee}\rangle\geq 0$ for all $\alpha \in S$. A dominant weight $\lambda \in
X({\bf P})$ is called $\bf P$-regular if $\langle \lambda ,\alpha
^{\vee}\rangle> 0$ for all $\alpha \notin I$, where ${\bf P} = {\bf P}_I$ is
a parabolic subgroup. A weight $\lambda$ defines a line bundle $\Ll
_{\lambda}$ on ${\bf G}/{\bf B}$. Line bundles on ${\bf G}/{\bf B}$
that correspond to dominant weights are ample. If a weight $\lambda$ is ${\bf
P}$-regular then the corresponding line bundle is ample on ${\bf G}/{\bf P}$.\\

Recall that the prime $p$ is {\it a good prime} for $\bf G$ if $p$ is
coprime to all the coefficients of the highest root of $\bf G$ written
in terms of the simple roots. In particular, if $\bf G$ is a simple
group of type $\bf A$ then all primes are good for $\bf G$; if
$\bf G$ is either of the type $\bf B$ or $\bf D$ then good primes are $p\geq 3$.\\

We will often use the Kempf vanishing theorem \cite{Kem}:

\begin{theorem}\label{th:Kempfvan}
Let $\bf P$ be a parabolic subgroup of $\bf G$, and $\Ll$ an effective
line bundle on ${\bf G}/{\bf P}$. Then ${\rm H}^{i}({\bf G}/{\bf P},\Ll)=0$ for $i>0$.
\end{theorem}

Next series of statements concerns line bundles on cotangent bundles.
We start with the vanishing theorem by Kumar et al. (\cite{KLT}, Theorem 5):

\begin{theorem}\label{th:KLTvantheorem}
Let ${\rm T}^{\ast}({\bf G}/{\bf P})$ be the total space of the
  cotangent bundle of a homogeneous space ${\bf G}/{\bf P}$, and
  $\pi \colon {\rm T}^{\ast}({\bf G}/{\bf P})\rightarrow {\bf G}/{\bf P}$ the projection.
Assume that $\mbox{\rm char} \ k$ is a good prime for $\bf G$. Let $\lambda\in
X({\bf P})$ be a $\bf P$-regular weight. Then 
\begin{equation}\label{eq:GmodPvanishing}
{\rm H}^i({\rm T}^*({\bf G}/{\bf P}), \pi^*\Ll (\lambda)) = 
{\rm H}^{i}({\bf G}/{\bf P},\Ll _{\lambda}\otimes \pi _{\ast}\Oo _{{\rm T}^{\ast}({\bf G}/{\bf P})}) = 0,
\end{equation}

for $i>0$.
\end{theorem}

Here $\Ll _{\lambda}$ denotes a line bundle that corresponds to the
weight $\lambda$. In particular, one has:
\begin{equation}\label{eq:OacyclicGmodB}
{\rm H}^i({\rm T}^*({\bf G}/{\bf B}),\Oo _{{\rm T}^*({\bf G}/{\bf B})}) = 0,\end{equation}

for $i>0$.

It was proved in \cite{Don} and \cite{MVdK} that nilpotent orbits in type
${\bf A}_n$ are normal and have rational singularities. This implies
(Propositions 4.6 and 4.9 of \cite{MVdK}):

\begin{theorem}\label{th:typeAvantheorem}
Let ${\bf G} = {\bf SL}_n(k)$, and ${\bf P}\subset {\bf G}$ a
parabolic subgroup. Then
\begin{equation}\label{eq:OacyclicGmodB}
{\rm H}^i({\rm T}^*({\bf G}/{\bf P}),\Oo _{{\rm T}^*({\bf G}/{\bf P})}) = 0,
\end{equation}

for $i>0$.
\end{theorem}

Yet another vanishing theorem was also proved in
\cite{KLT} (Theorem 6):

\begin{theorem}\label{th:KLTsubregnilpvarth}
Let ${\bf P}_{\alpha}\subset {\bf G}$ be a minimal parabolic subgroup corresponding to
a simple short root $\alpha$ in the root system of $\bf G$, and $p$ a good
prime for $\bf G$. Then 
\begin{equation}\label{eq:OacyclicGmodB}
{\rm H}^i({\rm T}^*({\bf G}/{\bf P}_{\alpha}),\Oo _{{\rm T}^*({\bf G}/{\bf P}_{\alpha})}) = 0,
\end{equation}

\noindent for $i>0$.
\end{theorem}

\begin{remark}
{\rm One believes that the higher cohomology of $\Oo _{{\rm T}^*({\bf
    G}/{\bf P})}$ must vanish for any parabolic subgroup $\bf P$ and
  suitable $p$ (\cite{BrKum}, p. 182). More generally, the vanishing
as in (\ref{eq:GmodPvanishing}) should hold for any dominant line
  bundle $\Ll _{\lambda}$.}
\end{remark}

\vspace{0.2cm}

Let now ${\bf G}$ be either of the type ${\bf B}$ or ${\bf D}$, and ${\bf
  P}$ a maximal parabolic subgroup of $\bf G$ such that the
  grassmannian ${\bf G}/{\bf P}$ is isomorphic to a smooth quadric. We will need a similar vanishing result as in
  Theorem \ref{th:KLTsubregnilpvarth} for such grassmannians.

\begin{lemma}\label{lem:quadricvanishing}
\quad Let ${\sf Q}_n$ be a smooth quadric of dimension $n$. Assume that $p$ 
is odd. Then ${\rm H}^i({\rm T}^*({\sf Q}_n),\Oo _{{\rm T}^*({\sf Q}_n)})$ =
$0$ for $i>0$.

\begin{proof}
Let $j:{\sf Q}_n\hookrightarrow \Pp (V)$ be the embedding of ${\sf Q}_n$ into the projective space $\Pp
(V)$, the dimension of $V$ being equal to $n+2$. Consider the adjunction sequence
\begin{equation}\label{eq:adjsequenceonquadric}
0\rightarrow \T _{{\sf Q}_n}\rightarrow j^{\ast}\T _{\Pp (V)}\rightarrow
\Oo _{{{\sf Q}_n}}(2)\rightarrow 0,
\end{equation}

and tensor it with $\Oo  _{{{\sf Q}_n}}(-1)$:

\begin{equation}\label{eq:adjseqonquadrictwisted-1}
0\rightarrow \T _{{\sf Q}_n}(-1)\rightarrow j^{\ast}\T _{\Pp (V)}(-1)\rightarrow
\Oo _{{{\sf Q}_n}}(1)\rightarrow 0.
\end{equation}

Note that $\Hom (j^{\ast}\T _{\Pp (V)}(-1),\Oo _{{{\sf
      Q}_n}}(1)) = {\rm H}^0({\sf Q}_n, j^{\ast}\Omega ^1 _{\Pp
      (V)}(2)) = {\rm Ker}(V^{\ast}\otimes
      V^{\ast}\rightarrow {\sf S}^2V^{\ast}/\langle q\rangle)$ (the
      last isomorphism comes from the Euler exact sequence on $\Pp (V)$
      restricted to ${\sf Q}_n$). Here $q\in {\sf S}^{2}{\sf
      V}^{\ast}$ is the quadratic form that defines the quadric ${\sf
      Q}_n$ (since $p$ is odd, we can identify quadratic and bilinear
      forms via polarization) and the map $j^{\ast}\T _{\Pp (V)}(-1)\rightarrow
\Oo _{{{\sf Q}_n}}(1)$ in (\ref{eq:adjseqonquadrictwisted-1}) corresponds to $q$. The kernel of this map, which is isomorphic
      to $\T _{{\sf Q}_n}(-1)$, is also equipped with a symmetric form
      that makes the bundle $\T _{{\sf Q}_n}(-1)$ an orthogonal vector
      bundle with trivial determinant. We then obtain an
      isomorphism of vector bundles $\T _{{\sf Q}_n}(-1)$ and $\Omega
      _{{\sf Q}_n}^1(1)$ (in fact, the fiber of $\T _{{\sf Q}_n}(-1)$
      over a point on ${\sf Q}_n$ is isomorphic to $l^{\perp}/l$,
      where $l$ is the line corresponding to this point and
      $l^{\perp}$ is the orthogonal complement to $l$ with respect to
      $q$). Therefore, $\T _{{\sf Q}_n}=\Omega
      _{{\sf Q}_n}^1(2)$. Dualizing the sequence
      (\ref{eq:adjsequenceonquadric}) and tensoring it then with $\Oo
      _{{\sf Q}_n}(2)$ we get:

\begin{equation}\label{eq:KosresfortangbuntoQ_n}
0\rightarrow \Oo _{{{\sf Q}_n}}\rightarrow j^{\ast}\Omega ^1_{\Pp (V)}(2)\rightarrow \T _{{\sf Q}_n}\rightarrow 0.
\end{equation}

Let $\pi : {\rm T}^*({\sf Q}_n)\rightarrow {\sf Q}_n$ be the 
projection. For any $i\geq 0$ one has an isomorphism ${\rm H}^{i}({\rm T}^*({\sf
  Q}_n),\Oo _{{\rm T}^*({\sf Q}_n)})$ = ${\rm H}^{i}({\sf Q}_n,{\sf
  S}^{\bullet}\T _{{\sf Q}_n})$, the projection $\pi$ being an affine
morphism. One has ${\rm H}^i({\sf Q}_n,\Oo _{{\sf Q}_n})=0$ for
$i>0$. Let $k\geq 1$. The Koszul resolution (Section
\ref{subsec:Kosres}) for the bundle ${\sf S}^k\T _{{\sf Q}_n}$ looks as follows:
\begin{equation}
0 \rightarrow j^{\ast}{\sf S}^{k-1}(\Omega ^1_{\Pp (V)}(2))\rightarrow j^{\ast}{\sf S}^{k}(\Omega ^1_{\Pp (V)}(2))\rightarrow {\sf S}^{k}\T _{{\sf Q}_n}\rightarrow 0.
\end{equation}

Consider the projective bundle $\Pp (\T _{\Pp (V)}(-1))$ over $\Pp
(V)$. It is isomorphic to the variety of partial flags ${\sf
  Fl}_{1,2,n+1}$. Let $\pi$ and $q$ be the projections of ${\sf
  Fl}_{1,2,n+1}$ onto $\Pp (V)$ and ${\sf Gr}_{2,n+1}$,
respectively. Note that for $k\geq 0$ one has ${\rm R}^{\bullet}\pi _{\ast}\Oo _{\pi}(k)={\sf S}^k\Omega
^1_{\Pp (V)}(k)$, where $\Oo _{\pi}(1)$ is the relatively ample
invertible sheaf on $\Pp (\T _{\Pp (V)}(-1))$. Thus, by the
projection formula, ${\sf S}^{k}(\Omega ^1_{\Pp (V)}(2)) = {\rm R}^{\bullet}\pi _{\ast}(\Oo _{\pi}(k)\otimes \pi ^{\ast}\Oo _{\Pp
  (V)}(k)) = {\rm R}^{\bullet}\pi
_{\ast}q^{\ast}\Oo _{{\sf Gr}_{2,n+1}}(k)$. \ On \ the \ other \ hand,
\ ${\rm H}^{i}(\Pp (V), {\rm R}^{\bullet}\pi
_{\ast}q^{\ast}\Oo _{{\sf Gr}_{2,n+1}}(k))$ \ = \ ${\rm H}^{i}(\Pp (\T
_{\Pp (V)}(-1),q^{\ast}\Oo _{{\sf Gr}_{2,n+1}}(k))$. Using the Kempf vanishing, we get ${\rm H}^{i}(\Pp (\T
_{\Pp (V)}(-1)),q^{\ast}\Oo _{{\sf Gr}_{2,n+1}}(k))$ = $0$ for $i>0$, the line bundle
$q^{\ast}\Oo _{{\sf Gr}_{2,n+1}}(k)$ being effective. Now tensor the short exact sequence
\begin{equation}\label{eq:quadbasicseq}
0\rightarrow \Oo _{\Pp (V)}(-2)\rightarrow \Oo _{\Pp (V)}\rightarrow
j_{\ast}\Oo _{{\sf Q}_n}\rightarrow 0,
\end{equation}

with ${\sf S}^{k}(\Omega ^1_{\Pp (V)}(2))$. The long exact cohomology
sequence gives:
\begin{eqnarray}\label{eq:lcohseqquadbasicseq}
& \dots \rightarrow {\rm H}^i(\Pp (V),{\sf S}^{k}(\Omega ^1_{\Pp
  (V)}(2))\otimes \Oo _{\Pp (V)}(-2))\rightarrow {\rm H}^i(\Pp
  (V),{\sf S}^{k}(\Omega ^1_{\Pp (V)}(2)))\rightarrow \\
& \rightarrow {\rm H}^i({\sf Q}_n,j^{\ast}{\sf S}^{k}(\Omega ^1_{\Pp (V)}(2))\rightarrow \dots . &\nonumber
\end{eqnarray}

By the above, the cohomology groups in the middle term vanish for $i>0$. As
 for the first term, one has an isomorphism ${\rm H}^i(\Pp (V),{\sf S}^{k}(\Omega ^1_{\Pp
  (V)}(2))\otimes \Oo _{\Pp (V)}(-2))={\rm H}^i({\sf Fl}_{1,2,n+1}, q^{\ast}\Oo_{{\sf
 Gr}_{2,n+1}}(k)\otimes \pi ^{\ast}\Oo _{\Pp (V)}(-2))$. The latter
 group is isomorphic to ${\rm H}^i({\sf Gr}_{2,n+1},\Oo_{{\sf
 Gr}_{2,n+1}}(k)\otimes {\rm R}^{\bullet}q_{\ast}\pi ^{\ast}\Oo _{\Pp
 (V)}(-2))$. It is easy to see (for example, using the Koszul resolution of the
structure sheaf of ${\sf Fl}_{1,2,n+1}$ inside the product $\Pp
(V)\times {\sf Gr}_{2,n+1}$) that ${\rm R}^{\bullet}q_{\ast}\pi
^{\ast}\Oo _{\Pp (V)}(-2) = \Oo _{{\sf
 Gr}_{2,n+1}}(-1)[-1]$. Thus, one obtains ${\rm H}^{i}({\sf
  Fl}_{1,2,n+1},q^{\ast}\Oo _{{\sf Gr}_{2,n+1}}(k)\otimes \pi
^{\ast}\Oo _{\Pp (V)}(-2)) = {\rm H}^{i-1}({\sf Gr}_{2,n+1},\Oo _{{\sf Gr}_{2,n+1}}
(k-1))=0$ for $i\neq 1$. From the sequence
 (\ref{eq:lcohseqquadbasicseq}) one has ${\rm H}^{i}({\sf Q}_n,j^{\ast}{\sf S}^{k}(\Omega ^1_{\Pp (V)}(2)) = 0$ for $i>0$. Finally, the sequence
 (\ref{eq:KosresfortangbuntoQ_n}) implies ${\rm H}^{i}({\sf Q}_n,{\sf
 S}^{k}\T _{{\sf Q}_n}) = 0$ for $i>0$, hence the statement.

\end{proof}
\end{lemma}


\section{\bf Cohomology of the Frobenius neighborhoods}


Let $X$ be a smooth variety over $k$. To compute the 
    cohomology groups ${\rm H}^{i}(X',{\mathcal End}({{\sf
    F}}_{\ast}\Oo _X))$ we will use the properties of sheaves
    ${\rm D}_X$ from Section \ref{subsec:Berthdiffoper}. Keeping
    the previous notation, we get:
\begin{eqnarray}\label{eq:chainofisom}
& {\rm H}^{j}(X',{\mathcal End}({\sf F}_{\ast}\Oo _X)) = 
{\rm H}^{j}(X',i^{\ast}\DD _X) = 
{\rm H}^{j}({\rm T}^{\ast}(X'),i_{\ast}i^{\ast}\DD _X) = & \nonumber \\
& = {\rm H}^{j}({\rm T}^{\ast}(X'),\DD _X\otimes i_{\ast}\Oo _{X^{'}}), 
\end{eqnarray}

the last isomorphism in (\ref{eq:chainofisom}) follows from
  the projection formula. Recall the projection $\pi \colon {\rm
    T}^{\ast}(X^{'})\rightarrow X^{'}$. Consider the bundle $\pi ^{\ast}\T
  ^{\ast}_{X^{'}}$. There is a tautological section $s$ of this bundle
such that the zero locus of $s$ coincides with $X^{'}$. Hence, one obtains the Koszul resolution:
\begin{equation}\label{eq:Kosres}
0\rightarrow \dots \rightarrow \bigwedge^{k}(\pi ^{\ast}\T _{X^{'}}) \rightarrow
\bigwedge^{k-1}(\pi ^{\ast}\T _{X^{'}})\rightarrow \dots \rightarrow \Oo
_{{\rm T}^{\ast}(X^{'})}\rightarrow i_{\ast}\Oo _{X^{'}}\rightarrow 0.
\end{equation}

Tensor the resolution (\ref{eq:Kosres}) with the sheaf
$\DD _X$. The rightmost cohomology group in (\ref{eq:chainofisom})
can be computed via the above Koszul resolution. 
\begin{lemma}\label{lem:cohofFrobneigh}
\ \ Fix \ \ $k\geq 0$. \ \ For \ any \ $j\geq 0$, \ if \ ${\rm H}^{j}({\rm
  T}^{\ast}(X),{\sf F}^{\ast}\bigwedge ^{k}(\pi ^{\ast}\T _{X'})) = 0$
  \ then \ ${\rm H}^{j}({\rm T}^{\ast}(X'),\DD _X\otimes \bigwedge ^{k}(\pi ^{\ast}\T _{X'})) = 0$. 
\end{lemma}

\begin{proof}
Denote $C^{k}$ the sheaf $\DD _X \otimes \bigwedge ^{k}(\pi ^{\ast}\T _{X'})$.
Take the direct image of $C^k$ with respect to $\pi$. Using the
projection formula we get:
\begin{eqnarray}
& {\rm R}^{\bullet}\pi _{\ast}C^{k} = {\rm
  R}^{\bullet}\pi _{\ast}(\DD _X\otimes \bigwedge ^{k}(\pi ^{\ast}\T
  _{X'})) = \pi _{\ast}(\DD _X\otimes \bigwedge ^{k}(\pi ^{\ast}\T
  _{X'})) & \nonumber \\
& = {\sf F}_{\ast}{\rm D} _X\otimes \bigwedge ^{k}(\T _{X'}), 
\end{eqnarray}

the morphism $\pi$ being affine. Hence,
\begin{equation}
{\rm H}^{j}({\rm T}^{\ast}(X'),\DD _X\otimes
\bigwedge ^{k}(\pi ^{\ast}\T _{X'})) = {\rm H}^{j}(X', {\sf
 F}_{\ast}{\rm D} _X\otimes \bigwedge ^{k}(\T _{X'})).
\end{equation}

The sheaf ${\sf F}_{\ast}{\rm D} _X\otimes
\bigwedge ^{k}(\T _{X'})$ is equipped with a filtration that is induced
  by the filtration on ${\sf F}_{\ast}{\rm D}_X$, the associated sheaf
  being isomorphic to ${\rm gr}({\sf F}_{\ast}{\rm D} _X)\otimes \bigwedge ^{k}(\T
_{X'}) = {\sf F}_{\ast}\pi _{\ast}\Oo _{{\rm T}^{\ast}(X)}\otimes
\bigwedge ^{k}(\T _{X'})$. Clearly, for $j\geq 0$ 
\begin{equation}
{\rm H}^{j}(X',{\rm gr}({\sf F}_{\ast}{\rm D} _X)\otimes
\bigwedge ^{k}(\T _{X'})) = 0 \quad \Rightarrow \quad {\rm H}^{j}(X',{\sf F}_{\ast}{\rm D}_{X}\otimes
\bigwedge ^{k}(\T _{X'})) = 0. 
\end{equation}

There are isomorphisms:
\begin{eqnarray}
& {\rm H}^{j}(X',{\sf F}_{\ast}\pi _{\ast}\Oo _{{\rm T}^{\ast}(X)}\otimes
\bigwedge ^{k}(\T _{X'})) = {\rm H}^{j}(X,\pi _{\ast}\Oo _{{\rm
  T}^{\ast}(X)}\otimes {\sf F}^{\ast}\bigwedge ^{k}(\T _{X'})) = & \nonumber \\
& = {\rm H}^{j}({\rm T}^{\ast}(X),\pi ^{\ast}{\sf F}^{\ast}\bigwedge ^{k}(\T
_{X'})),
\end{eqnarray}

the \ last \ isomorphism \ following  \ from \ the \ projection \
formula. \ Finally, \ ${\rm H}^{j}({\rm T}^{\ast}(X),\pi ^{\ast}{\sf F}^{\ast}\bigwedge ^{k}(\T
_{X'})) = {\rm H}^{j}({\rm T}^{\ast}(X),{\sf F}^{\ast}\bigwedge ^{k}(\pi
^{\ast}\T _{X'}))$, hence the statement of the lemma. 
\end{proof}

\begin{remark}
Assume that for a given $X$ one has ${\rm H}^{j}({\rm T}^{\ast}(X),{\sf F}^{\ast}\bigwedge
^{k}(\pi ^{\ast}\T _{X'})) = 0$ for $j>k\geq 0$. The \ spectral \ sequence \
${\rm E}_1^{p,q} = {\rm H}^{p}({\rm T}^{\ast}(X'),\DD _X\otimes
\bigwedge ^q(\pi ^{\ast}\T _{X'}))$ \ converges \ to \ ${\rm H}^{p-q}({\rm
  T}^{\ast}(X'),\DD _X\otimes i_{\ast}\Oo _{X^{'}})$.  Lemma \ref{lem:cohofFrobneigh} and the resolution
(\ref{eq:Kosres}) then imply:

\begin{equation}
{\rm H}^{j}({\rm
  T}^{\ast}(X'),\DD _X\otimes i_{\ast}\Oo _{X^{'}}) = {\rm H}^{j}(X',{\mathcal End}({\sf F}_{\ast}\Oo _X)) = 0
\end{equation}

for $j>0$, and 

\begin{equation}
{\rm H}^{j}({\rm  T}^{\ast}(X),{\sf F}^{\ast}i_{\ast}\Oo _{X'}) = 0
\end{equation}
 for $j>0$.
\end{remark}

\begin{remark}
{\rm  The complex ${\tilde C}^{k}\colon = {\sf F}^{\ast}\pi ^{\ast}\bigwedge^{k}(\T
_{X'})$ is quasiisomorphic to ${\sf F}^{\ast}i_{\ast}\Oo _{X'}$, the
structure sheaf of the Frobenius neighborhood of the zero
section $X'$ in the cotangent bundle ${\rm T}(X')$. Below we show that
if $X$ is either a smooth quadric (for odd $p$) or a
partial flag variety then 
\begin{equation}
{\rm H}^{j}({\rm T}^{\ast}(X),{\sf F}^{\ast}i_{\ast}\Oo _{X'}) = 0
\end{equation}

for $j > 0$}. 
\end{remark}


\section{\bf Quadrics}\label{sec:3}


\begin{theorem}\label{th:quadrics}
Let ${\sf Q}_n$ be a smooth quadric of dimension $n$. \ Assume \ that \ $p$ \ is \
odd. \ Then \ ${\rm H}^{i}({\sf Q}_n,{\mathcal End}_{\Oo ^{p}_{{\sf Q}_n}}(\Oo _{{\sf Q}_n})) = 0 $ for $i > 0$.
\end{theorem}

\begin{proof}
Let $V$ be a vector space of dimension $n+2$. A smooth quadric ${\sf Q}_n\subset
\Pp (V)$ is a homogeneous space, and a hypersurface of degree two in
$\Pp (V)$. To keep the notation simple, we will ignore the Frobenius twist of
${\sf Q}_n$ and will be dealing with the absolute Frobenius morphism
${\sf F}:{\sf Q}_n\rightarrow {\sf Q}_n$; this will not change the
cohomology groups in question.

Lemma \ref{lem:auxlemmaquadr} shows that
${\rm H}^{i}({{\sf Q}_n},\bigwedge^{r}\T _{{{\sf Q}_n}}\otimes  {\sf F}_{\ast}\pi _{\ast}\Oo _{{\rm
  T}^{\ast}({{\sf Q}_n})}) = 0$ for $i>r\geq 0$. Lemma
\ref{lem:cohofFrobneigh} then gives
${\rm H}^{i}({\rm T}^{\ast}({{\sf Q}_n}),\DD _{{\sf Q}_n}\otimes
\bigwedge^{r}(\pi ^{\ast}\T _{{{\sf Q}_n}})) = 0$ for $i>r\geq
0$. Using Remark 2, we get:
\begin{equation}
{\rm H}^i({{\sf Q}_n},{\mathcal End}_{\Oo ^{p}_{{\sf Q}_n}}(\Oo _{{\sf Q}_n})) = 0
\end{equation}

for $i>0$. 
\end{proof}

\begin{lemma}\label{lem:auxlemmaquadr}
Let $\Ll$ be an effective line bundle on ${\sf Q}_n$. Then for $i>r\geq 0$ one has 
${\rm H}^{i}({{\sf Q}_n},\bigwedge^{r}\T _{{{\sf Q}_n}}\otimes {\sf F}_{\ast}\pi _{\ast}\Oo _{{\rm  T}^{\ast}({{\sf Q}_n})}\otimes \Ll) = 0$.
\end{lemma}

\begin{proof}
Recall that $\pi _{\ast}\Oo _{{\rm  T}^{\ast}({{\sf Q}_n})}={\sf
  S}^{\bullet}\T _{{\sf Q}_n}$. To simplify the notation, put also ${\mathbb L} = {\sf F}_{\ast}{\sf S}^{\bullet}\T _{{{\sf Q}_n}}\otimes \Ll$.
If $r=0$ and $\Ll$ is ample then the statement follows from Theorem
\ref{th:KLTvantheorem}. Indeed, by the projection formula one has:
\begin{equation}\label{eq:r=0k>0quad}
{\rm H}^{i}({{\sf Q}_n},{\mathbb L}) = {\rm H}^{i}({{\sf Q}_n},{\sf S}^{\bullet}\T
_{{\sf Q}_n}\otimes \Ll ^{\otimes p}) =0
\end{equation}

for $i>0$. If $\Ll =\Oo _{{\sf Q}_n}$ then Lemma \ref{lem:quadricvanishing} implies:
\begin{equation}\label{eq:Ovanishforquadrics}
{\rm H}^{i}({{\sf Q}_n},{\mathbb L}) = {\rm H}^{i}({{\sf Q}_n},{\sf S}^{\bullet}\T _{{\sf Q}_n}) = 0,
\end{equation}

for $i>0$. Let $r\geq 1$. We argue by
induction on $r$, the base of induction being $r = 1$. Denote $j:{\sf Q}_n\hookrightarrow \Pp (V)$
the embedding. Recall the adjunction sequence (see Lemma \ref{lem:quadricvanishing})
\begin{equation}\label{eq:adjseqonquadric}
0\rightarrow \T _{{\sf Q}_n}\rightarrow j^{\ast}\T _{\Pp (V)}\rightarrow
\Oo _{{{\sf Q}_n}}(2)\rightarrow 0,
\end{equation}

and consider the Euler sequence on $\Pp (V)$ restricted to ${{\sf Q}_n}$:
\begin{equation}\label{eq:Eulerseqrestonquad}
0\rightarrow \Oo _{{\sf Q}_n}\rightarrow V\otimes \Oo _{{\sf Q}_n}(1)\rightarrow j^{\ast}\T _{\Pp (V)}\rightarrow 0.
\end{equation}

The sequence (\ref{eq:adjseqonquadric}) gives rise to a short exact
sequence (see Section \ref{sec:prelim}):
\begin{equation}\label{eq:extpowerEulerquadr}
0\rightarrow \bigwedge^{r}\T _{{\sf Q}_n}\rightarrow j^{\ast}\bigwedge^{r}{\T  _{\Pp (V)}}\rightarrow \bigwedge^{r-1}\T _{{\sf Q}_n}\otimes \Oo _{{\sf Q}_n}(2)\rightarrow 0.
\end{equation}

Consider first the case $r=1$. Tensoring the sequence (\ref{eq:adjseqonquadric}) with ${\mathbb L}$, one obtains:
\begin{equation}\label{eq:tensoradjwith}
0\rightarrow \T _{{{\sf Q}_n}}\otimes {\mathbb L} \rightarrow
\T _{\Pp (V)}\otimes {\mathbb L}\rightarrow \Oo _{{{\sf Q}_n}}(2)\otimes {\mathbb L}\rightarrow 0.
\end{equation}

As in (\ref{eq:r=0k>0quad}), one has ${\rm H}^i({\sf Q}_n,\Oo _{{{\sf
      Q}_n}}(2)\otimes {\mathbb L})=0$ for $i>0$. Tensoring the sequence (\ref{eq:Eulerseqrestonquad}) with ${\mathbb L}$,
  we get:
\begin{equation}\label{eq:EulerseqtenswithSTandL}
 0\rightarrow {\mathbb L}\rightarrow
  V\otimes \Oo _{{{\sf Q}_n}}(1)\otimes {\mathbb L}\rightarrow \T _{\Pp (V)}\otimes {\mathbb L}\rightarrow 0. 
\end{equation}

Again by (\ref{eq:r=0k>0quad}) and (\ref{eq:Ovanishforquadrics}), one
has ${\rm H}^i({\sf Q}_n,{\mathbb L})={\rm H}^i({\sf Q}_n, V\otimes
\Oo _{{{\sf Q}_n}}(1)\otimes {\mathbb L})=0$ for $i>0$. Hence, ${\rm H}^{i}({{\sf Q}_n},\T _{\Pp (V)}\otimes {\mathbb L}) =
  0$ for $i > 0$. From (\ref{eq:tensoradjwith}) we conclude that 
${\rm H}^{i}({{\sf Q}_n},\T _{{{\sf Q}_n}}\otimes {\mathbb L}) = 0$
  for $i > 1$.

Now fix $m>1$. Assume that for $r \leq m$ one has 
${\rm H}^{i}({{\sf Q}_n},\bigwedge^{r}\T _{{{\sf Q}_n}}\otimes {\mathbb L}) = 0$ for $i
  > r$. Let us prove that ${\rm H}^{i}({{\sf Q}_n},\bigwedge^{r+1}\T
  _{{{\sf Q}_n}}\otimes {\mathbb L}) = 0$ for $i >
r+1$. Tensoring the sequence (\ref{eq:extpowerEulerquadr}) (for $r+1$) with ${\mathbb L}$, we get:
\begin{equation}\label{eq:extpowerEulerquadrtensored}
0\rightarrow \bigwedge^{r+1}\T _{{{\sf Q}_n}}\otimes {\mathbb L}\rightarrow
\bigwedge^{r+1}\T _{\Pp (V)}\otimes {\mathbb L}\rightarrow
\bigwedge^{r}\T _{{{\sf Q}_n}}\otimes \Oo _{{{\sf Q}_n}}(2)\otimes {\mathbb L}\rightarrow 0.
\end{equation}

By the inductive assumption one has ${\rm H}^{i}({{\sf Q}_n},\bigwedge^{r}\T _{{{\sf Q}_n}}\otimes \Oo
  _{{{\sf Q}_n}}(2)\otimes {\mathbb L}) = 0$ for $i > r$. Let us show that ${\rm
  H}^{i}({{\sf Q}_n},\bigwedge^{r+1}\T _{\Pp (V)}\otimes {\mathbb L}) = 0$ for $i > 0$. Clearly, this will imply ${\rm H}^{i}({{\sf Q}_n},\bigwedge^{r+1}\T
  _{{{\sf Q}_n}}\otimes {\mathbb L}) = 0$ for $i > r+1$. From the
  Euler sequence on $\Pp (V)$ one obtains:
\begin{equation}
0\rightarrow \bigwedge^{l-1}\T _{\Pp (V)}\rightarrow
  \bigwedge^{l}V\otimes \Oo _{\Pp (V)}(l)\rightarrow \bigwedge^{l}\T _{\Pp (V)}\rightarrow 0.
\end{equation}

Gluing together these short exact sequences for $l = 1,\dots r+1$,
one arrives at the long exact sequence:
\begin{equation}\label{eq:kthextpowerofEulerseqrestonquad}
0\rightarrow \Oo _{\Pp (V)}\rightarrow V\otimes \Oo _{\Pp
  (V)}(1)\rightarrow \bigwedge^{2}V\otimes \Oo _{\Pp (V)}(2)\rightarrow
  \dots \rightarrow \bigwedge^{r+1}\T _{\Pp (V)}\rightarrow 0.
\end{equation}

Tensoring the above sequence with ${\mathbb L}$, and using Theorem \ref{th:KLTvantheorem} and
  Lemma \ref{lem:quadricvanishing}, we conclude that ${\rm
  H}^{i}({{\sf Q}_n},\bigwedge^{r+1}\T _{\Pp (V)}\otimes {\mathbb L})
  = 0$ for $i > 0$. Hence the statement.
\end{proof}


\section{\bf Partial flag varieties}\label{sec:4}


Let $V$ be a vector space over $k$ of dimension $n+1$, and $X$ 
the flag variety ${\rm F}_{1,n,n+1}$, a smooth divisor in $\Pp
(V)\times \Pp (V^{\ast})$ of bidegree $(1,1)$. As in the previous
section, we omit the superscript at $X$ meaning the Frobenius twist
and consider the absolute Frobenius morphism.
\begin{theorem}\label{th:partialflags}
${\rm H}^{i}(X,{\mathcal End}_{\Oo ^{p}_X}(\Oo _X)) = 0 $ for $i > 0$.
\end{theorem}

\begin{proof} 
The proof is similar to that of Theorem \ref{th:quadrics}.
Consider the line bundle $\Oo _{\Pp (V)}(1)\boxtimes \Oo
  _{\Pp (V^{\ast})}(1)$ over $\Pp (V)\times \Pp
  (V^{\ast})$ (the symbol $\boxtimes$ denotes the external tensor
  product). Then $X$ is isomorphic to the zero locus of a section 
of  $\Oo _{\Pp (V)}(1)\boxtimes \Oo _{\Pp (V^{\ast})}(1)$. Consider
the adjunction sequence:
\begin{equation}\label{eq:adjsequenceincquad}
0\rightarrow \T _X\rightarrow \T _{\Pp (V)\times \Pp
  (V^{\ast})}\otimes \Oo _X\rightarrow \Oo _{X}(1)\boxtimes \Oo
  _{X}(1)\rightarrow 0.
\end{equation}

\vspace*{0.2cm}

Note that $\T _{\Pp (V)\times \Pp (V^{\ast})} = \T _{\Pp (V)}\boxplus
\T _{\Pp (V^{\ast})}$, the symbol $\boxplus$ denoting the external direct
sum. One has the Euler sequences:
\begin{equation}\label{eq:Eulerseq1}
0\rightarrow \Oo _{\Pp (V)}\rightarrow V\otimes \Oo _{\Pp
  (V)}(1)\rightarrow \T _{\Pp (V)}\rightarrow 0,
\end{equation}

 and 
\begin{equation}\label{eq:Eulerseq2}
0\rightarrow \Oo _{\Pp (V^{\ast})}\rightarrow V^{\ast}\otimes \Oo _{\Pp
  (V^{\ast})}(1)\rightarrow \T _{\Pp (V^{\ast})}\rightarrow 0.
\end{equation}

Hence the short exact sequence:
\begin{eqnarray}\label{eq:tangbunseqonincquad}
& 0\rightarrow \Oo _{\Pp (V)\times \Pp (V^{\ast})}\oplus \Oo _{\Pp (V)\times \Pp (V^{\ast})}\rightarrow
  V^{\ast}\otimes \Oo _{\Pp
  (V^{\ast})}(1)\boxplus \Oo _{\Pp (V)}(1)\otimes V\rightarrow & \nonumber \\
& \rightarrow \T _{\Pp (V)\times \Pp (V^{\ast})}\rightarrow 0 .
\end{eqnarray}

Show that ${\rm H}^{i}(X,\bigwedge^{r}\T _{X}\otimes  {\sf F}_{\ast}\pi _{\ast}\Oo
_{{\rm T}^{\ast}(X)}) = 0$ for $i > r\geq 0$. If $r=0$ then by Theorem \ref{th:typeAvantheorem} one has:
\begin{equation}\label{eq:typeAvanishforO}
{\rm H}^{i}(X,{\sf F}_{\ast}\pi _{\ast}\Oo _{{\rm
 T}^{\ast}(X)}) = {\rm H}^{i}(X,\pi _{\ast}\Oo
 _{{\rm T}^{\ast}(X)}) = {\rm H}^{i}(X,{\sf
 S}^{\bullet}\T _X) = 0
\end{equation}

for $i > 0$. Let now $r\geq 1$. The sequence (\ref{eq:adjsequenceincquad}) gives rise to a short exact
sequence:
\begin{equation}\label{eq:extpoweradjseqonincquad}
0\rightarrow \bigwedge^{r}\T _{X}\rightarrow \bigwedge^{r}{\T
  _{{\Pp (V)}\times {\Pp (V)}}}\otimes \Oo _{X}\rightarrow \wedge
  ^{r-1}\T _{X}\otimes (\Oo _{X}(1)\boxtimes \Oo _{X}(1))\rightarrow 0.
\end{equation}

We argue again by induction on $r$. Let first $r=1$.
Tensoring the sequence (\ref{eq:adjsequenceincquad}) with ${\sf F}_{\ast}{\sf S}^{\bullet}\T _{X}$, we get:
\begin{eqnarray}\label{eq:tensoradjwithII}
& 0\rightarrow \T _{X}\otimes {\sf F}_{\ast}{\sf S}^{\bullet}\T _{X}\rightarrow \T _{\Pp (V)\times \Pp
  (V^{\ast})}\otimes {\sf F}_{\ast}{\sf S}^{\bullet}\T _{X}\rightarrow & \nonumber \\
& \rightarrow (\Oo _{X}(1)\boxtimes \Oo _{X}(1))\otimes {\sf F}_{\ast}{\sf S}^{\bullet}\T _{X}\rightarrow 0.
\end{eqnarray}

By Theorem \ref{th:KLTvantheorem}, one has for $i > 0$:
\begin{equation}
{\rm H}^{i}(X,(\Oo _{X}(1)\boxtimes \Oo
  _{X}(1))\otimes {\sf F}_{\ast}{\sf S}^{\bullet}\T _{X}) = 
{\rm H}^{i}(X,(\Oo _{X}(p)\boxtimes \Oo _{X}(p))\otimes {\sf S}^{\bullet}\T _{X}) = 0.
\end{equation}

On the other hand, tensoring sequence
(\ref{eq:tangbunseqonincquad}) (considered first on ${\Pp (V)}\times
  {\Pp (V^{\ast}})$, and
  then restricted to $X$) with ${\sf F}_{\ast}{\sf
  S}^{\bullet}\T _{X}$, one obtains:
\begin{eqnarray}\label{eq:tangbunseqonincquadII}
& 0\rightarrow {\sf F}_{\ast}{\sf S}^{\bullet}\T _{X}^{\oplus 2}\rightarrow
  (V^{\ast}\otimes \Oo _{\Pp (V^{\ast})}(1)\boxplus
  \Oo _{\Pp (V)}(1)\otimes V)\otimes {\sf F}_{\ast}{\sf S}^{\bullet}\T
  _{X}\rightarrow & \nonumber \\
& \rightarrow \T _{\Pp (V)\times \Pp (V^{\ast})}\otimes {\sf F}_{\ast}{\sf S}^{\bullet}\T _{X}\rightarrow 0.
\end{eqnarray}

The leftmost term in the above sequence has vanishing higher
cohomology by (\ref{eq:typeAvanishforO}). The middle term is the
direct sum of several copies of the bundle ${\sf F}_{\ast}{\sf
  S}^{\bullet}\T _{X}$ tensored with an effective line bundle on $X$. Let $\Ll$ be an arbitrary effective line bundle
 on $X$ (i.e. the one isomorphic to either $\Oo _X(k)\boxtimes \Oo _X$
 or $\Oo _X\boxtimes \Oo _X(k)$ for some $k\geq 0$). For any $l\geq 1$ there is a filtration on ${\sf S}^{l}
(\T _{\Pp (V)\times \Pp (V^{\ast})}\otimes \Oo _X)$ that comes from
the sequence (\ref {eq:adjsequenceincquad}), the graded factors of this
 filtration being isomorphic to ${\sf S}^{i}\T _{X}\otimes (\Oo
 _{X}(l-i)\boxtimes \Oo _{X}(l-i))$ for $0\leq i\leq l$. Tensoring ${\sf S}^{l}
(\T _{\Pp (V)\times \Pp (V^{\ast})}\otimes \Oo _X)$ with $\Ll$, the graded factors ${\sf S}^{i}\T _{X}\otimes (\Oo
 _{X}(l-i)\boxtimes \Oo _{X}(l-i))$ are tensored with $\Ll$. For $i < l$ the
 higher cohomology of the corresponding graded factor vanish by
 Theorem \ref{th:KLTvantheorem}. On the other hand, the higher
 cohomology of ${\sf S}^{l}(\T _{\Pp (V)\times \Pp (V^{\ast})}\otimes
 \Oo _X)\otimes \Ll$ are easily seen to vanish (use the Koszul resolutions associated to the Euler sequences
 (\ref{eq:Eulerseq1}) and (\ref{eq:Eulerseq2}), and the Kempf
 vanishing). One obtains ${\rm H}^i(X,{\sf S}^{\bullet}\T _X\otimes \Ll) = 0$ for $i >
 1$. Hence, the cohomology of the middle term in
 (\ref{eq:tangbunseqonincquadII}) vanish for $i>1$, and, therefore, ${\rm H}^{i}(X,\T
 _{\Pp (V)\times \Pp (V^{\ast})}\otimes {\sf F}_{\ast}{\sf
   S}^{\bullet}\T _{X}) = 0$ for $i>1$. Coming back to the sequence
 (\ref{eq:tensoradjwithII}), we get ${\rm H}^{i}(X,\T _{X}\otimes
 {\sf F}_{\ast}{\sf S}^{\bullet}\T _{X}) = 0$ for $i>1$. The inductive
 step is completely analogous to that in the proof of Lemma
 \ref{lem:auxlemmaquadr} (one uses the Euler sequences
 (\ref{eq:Eulerseq1}) and (\ref{eq:Eulerseq2}) to ensure that ${\rm H}^{i}(X,\bigwedge^{r+1}\T _{\Pp (V)\times \Pp
   (V^{\ast})}\otimes {\sf F}_{\ast}{\sf
   S}^{\bullet}\T _{X}) = 0$ for $i > r+1$). This allows to complete the
 induction argument. 
\end{proof}

\section{\bf Derived equivalences}\label{sec:derequiv}

\begin{definition}\label{def:tiltingdefinition}
A coherent sheaf $\Ee$ on a smooth variety $X$ over an algebraically
closed field $k$ is called a {\sf tilting generator} of
the bounded derived category $\D ^b(X)$ of coherent sheaves on
$X$ if the following holds:
\begin{enumerate}
\item The sheaf $\Ee$ is a tilting object in $\D^b(X)$, that
  is, for any $i \geq 1$ one has $\Ext^i(\Ee,\Ee)=0$.
\item The sheaf $\Ee$ generates the derived category $\D ^{-}(X)$
  of complexes bounded from above, that is, if for some object $\Ff
  \in \D^{-}(X)$ one has $\RHom ^{\bullet}(\Ff ,\Ee)=0$, then $\Ff=0$.
\end{enumerate}
\end{definition}

Tilting generators are a tool to construct derived equivalences. One has:

\begin{lemma}$\rm ($\cite{Kaledquant}, Lemma 1.2$\rm )$\label{lem:tiltinglemma}
Let $X$ be a smooth scheme, $\Ee$ a tilting generator of the derived category $\D^b(X)$, and
denote $R={\mathcal End}(\Ee)$. Then the
algebra $R$ is left-Noetherian, and the correspondence $\Ff \mapsto
\RHom^{\bullet}(\Ee,\Ff)$ extends to an equivalence
\begin{equation}\label{equi}
  \D^b(X) \rightarrow \D^b(R-\mbox{\rm mod}^{\rm fg})
\end{equation}
between the bounded derived category $\D^b(X)$ of coherent
sheaves on $X$ and the bounded derived category $\D^b(R-\mbox{\rm mod}^{\rm fg})$ of
finitely generated left $R$-modules.
\end{lemma}


We now recall the derived localization theorem of \cite{BMR}. Let $\bf G$ be a semisimple algebraic group over $k$, ${\bf G}/{\bf
  B}$ the flag variety, and $\Uu ({\mathfrak g})$ the universal
  enveloping algebra of the corresponding Lie algebra. The center of
  $\Uu ({\mathfrak g})$ contains the ``Harish--Chandra part''
  ${\mathfrak Z}_{\rm HC}=\Uu ({\mathfrak g})^{\bf G}$. Denote $\Uu
  ({\mathfrak g})_0$ the central reduction $\Uu ({\mathfrak g})\otimes _{{\mathfrak
  Z}_{\rm HC}} {\bf k}$, where $\bf k$ is the trivial $\mathfrak g$-module.
Consider the category ${\rm D}_{{\bf G}/{\bf
B}}$-\mbox{\rm mod} of coherent ${\rm D}_{{\bf G}/{\bf B}}$-modules and the
category $\Uu ({\mathfrak g})_0$-\mbox{\rm mod} of finitely generated modules over $\Uu
  ({\mathfrak g})_0$. The
  derived localization theorem (Theorem 3.2, \cite{BMR}) states:

\begin{theorem}\label{th:derBB-equivalence}
Let char $k = p>h$, where $h$ is the Coxeter number of the group $\bf
G$. Then there is an equivalence of derived categories:
\begin{equation}\label{eq:derBB-equivalence}
\Dd ^{b}({\rm D}_{{\bf G}/{\bf B}}-\mbox{\rm mod}) \simeq \Dd ^{b}(\Uu ({\mathfrak g})_0-\mbox{\rm mod}),
\end{equation}

\end{theorem}

There are also ``unbounded'' versions of Theorem
\ref{th:derBB-equivalence}, the bounded categories in
(\ref{eq:derBB-equivalence}) being replaced by categories unbounded from
above or below (see Remark 2 on p. 18 of \cite{BMR}). Thus extended, Theorem \ref{th:derBB-equivalence} implies:

\begin{lemma}\label{lem:Frobgeneratorforp>h}
Let $\bf G$ be a semisimple algebraic group over $k$, and $X={\bf
  G}/{\bf B}$ the flag variety. Let char $k = p >h$, where $h$ is the
  Coxeter number of the group ${\bf G}$. Then the bundle ${\sf F}_{\ast}\Oo _X$
satisfies the condition (2) of Definition \ref{def:tiltingdefinition}.
\begin{proof}
We \ need \ to \ show \ that \ for \ an \ object \ $\Ee
  \in \D^{-}(X')$ \ the \ equality \ $\RHom ^{\bullet}(\Ee ,{\sf F}_{\ast}\Oo
  _X)=0$ implies $\Ee=0$. By adjunction we get:
\begin{equation}
{\mathbb H}^{\bullet}(X,({\sf F}^{\ast}\Ee)^{\vee}) = 0.
\end{equation}

On the other hand,
\begin{eqnarray}
& ({\sf F}^{\ast}\Ee)^{\vee}={\mathcal R}{\mathcal
  Hom}_{\Oo _X}({\sf F}^{\ast}\Ee,\Oo _X) = {\mathcal R}{\mathcal
  Hom}_{\Oo _X}({\sf F}^{\ast}\Ee,{\sf F}^{\ast}\Oo _{X'})  = & \nonumber \\
& ={\sf F}^{\ast}{\mathcal R}{\mathcal Hom}_{\Oo _{X'}}(\Ee,\Oo _{X'})={\sf F}^{\ast}\Ee ^{\vee}.
\end{eqnarray}

By the Cartier descent (see Subsection \ref{subsec:Cartier}), the object ${\sf F}^{\ast}\Ee ^{\vee}$ is an object of the category
$\Dd ^{-}(\Dd _X-\mbox{\rm mod})$. Now ${\sf F}^{\ast}\Ee ^{\vee}$ is annihilated by the functor
${\rm R}\Gamma$. Under our assumptions on $p$, this functor is an equivalence of categories by Theorem
\ref{th:derBB-equivalence}. Hence, ${\sf F}^{\ast}\Ee ^{\vee}$ is
quasiisomorphic to zero, and therefore so are $\Ee ^{\vee}$ and $\Ee$. 
\end{proof}
\end{lemma}

\begin{corollary}\label{cor:Frobgeneratorforp>h}
Let $\bf P$ be a parabolic subgroup of $\bf G$, and let $p>h$. Then
  the bundle ${\sf F}_{\ast}\Oo _{{\bf G}/{\bf P}}$ is a generator in ${\rm D}^{b}({\bf
  G}/{\bf P})$.
\end{corollary}

\begin{proof}
Denote $Y={\bf G}/{\bf P}$, and let $\pi : X = {\bf G}/{\bf B}\rightarrow Y$ be the projection. As before, one has to show
that for any object $\Ee \in \Dd ^{-}(Y')$ the equality $\RHom ^{\bullet}(\Ee ,{\sf F}_{\ast}\Oo
  _Y)=0$ implies $\Ee=0$. Notice that  ${\rm R}^{\bullet}\pi _{\ast}\Oo _{X} = \Oo _{Y}$. The condition $\RHom ^{\bullet}(\Ee ,{\sf F}_{\ast}\Oo
  _Y)=0$ then translates into:
\begin{equation}
\RHom ^{\bullet}(\Ee ,{\sf F}_{\ast}\Oo  _Y)= \RHom ^{\bullet}(\Ee,
{\sf F}_{\ast}{\rm R}^{\bullet}\pi _{\ast}\Oo _X) = \RHom ^{\bullet}(\Ee
,{\rm R}^{\bullet}\pi _{\ast}{\sf F}_{\ast}\Oo _X) =0.
\end{equation}

By adjunction we get:
\begin{equation}
\RHom ^{\bullet}(\pi ^{\ast}\Ee ,{\sf F}_{\ast}\Oo _X) =0.
\end{equation}

Lemma \ref{lem:Frobgeneratorforp>h} then implies that the
object $\pi ^{\ast}\Ee = 0$ in $\Dd ^{-}(X')$. Applying to it the
functor $\pi _{\ast}$ and using the projection formula, we get ${\rm
  R}^{\bullet}\pi _{\ast}\pi ^{\ast}\Ee=\Ee \otimes {\rm
  R}^{\bullet}\pi _{\ast}\Oo _{X'}=\Ee$, hence $\Ee = 0$, q.e.d.
\end{proof}

\begin{corollary}
Let $X$ be the variety of partial flags ${\rm F}_{1,n,n+1}$. Assume that $p>n+1$. Then ${\sf
  F}_{\ast}\Oo _X$ is a tilting bundle.
\end{corollary}

\begin{proof}
Follows from Theorem \ref{th:partialflags} and Corollary \ref{cor:Frobgeneratorforp>h}.
\end{proof}

\begin{corollary}
Let ${\sf Q}_n$ be a smooth quadric of dimension $n$. Assume that
  $p>n-(-1)^k$, where $k\equiv n(2)$. Then ${\sf
  F}_{\ast}\Oo _{{\sf Q}_n}$ is a tilting bundle.
\end{corollary}

\begin{proof}
Follows from Theorem \ref{th:quadrics} and Corollary \ref{cor:Frobgeneratorforp>h}.
\end{proof}

\begin{remark}
{\rm The derived categories of coherent sheaves on quadrics and partial
flag varieties ${\rm F}_{1,n,n+1}$ were described in \cite{Kap} in
terms of exceptional collections. By \cite{Lan}, decomposition of
the Frobenius pushforward of the structure sheaf on a quadric gives
rise to the full exceptional collection from \cite{Kap}. The relation
of the Frobenius pushforward of the structure sheaf on ${\rm
  F}_{1,n,n+1}$ and the full exceptional collection given in {\it
  loc.cit.} will be a subject of a future paper (for $n=2$ the variety ${\rm
  F}_{1,2,3}$ is isomorphic to the flag variety ${\bf SL}_3/{\bf B}$
and the decomposition of ${\sf F}_{\ast}\Oo _{{\bf SL}_3/{\bf B}}$
into the direct sum of indecomposables was
found in \cite{HKR})}. 
\end{remark}


\end{document}